\documentclass[a4paper,11pt]{article}
\usepackage[all,cmtip]{xy}
\newdir{ >}{{}*!/-5pt/@{>}}
\usepackage{amssymb,latexsym,amsmath}
\usepackage{amsthm}
\usepackage[a4paper]{geometry}
\usepackage[colorlinks=true,breaklinks=true,linkcolor=black,citecolor=black,hidelinks]{hyperref}
\usepackage{verbatim}
\usepackage[textsize=footnotesize,backgroundcolor=white]{todonotes}
\usepackage{times}
\usepackage{manfnt}
\usepackage{mathrsfs}
\usepackage{bm}
\usepackage{upgreek}

\theoremstyle{plain}

\swapnumbers
\theoremstyle{plain}
\newtheorem{theorem}{Theorem}[section]

\newtheorem{lemma}[theorem]{Lemma}
\newtheorem{prop}[theorem]{Proposition}
\newtheorem{corollary}[theorem]{Corollary}

\theoremstyle{definition}
\newtheorem{ntheorem}[theorem]{\theoremname}
\newcommand{\theoremname}{Test}
\newenvironment{des}[1]{\renewcommand{\theoremname}{#1} \begin{ntheorem}}{\end{ntheorem}}
\newcommand{\bdes}{\begin{des}}
\newcommand{\edes}{\end{des}}
\newtheorem{ex}[theorem]{Example}
\newtheorem{definition}[theorem]{Definition}
\newtheorem{exer}{Exercise}
\newtheorem{rem}[theorem]{Remark}
\newtheorem{observation}[theorem]{Observation}
\newtheorem{app}[theorem]{Application}
\newtheorem{his}[theorem]{Historical Note}
\newtheorem*{notation}{Notation}
\newtheorem*{conv}{Convention}
\newtheorem{Observation}[theorem]{Observation}
\newtheorem{cons}[theorem]{Construction}
\newtheorem{diss}[theorem]{}
\renewcommand{\theenumi}{(\roman{enumi})}

\newcommand{\bthm}{\begin{theorem}}
\newcommand{\ethm}{\end{theorem}}

\newcommand{\bprop}{\begin{prop}}
\newcommand{\eprop}{\end{prop}}

\newcommand{\bcoro}{\begin{corollary}}
\newcommand{\ecoro}{\end{corollary}}

\newcommand{\bex}{\begin{ex}}
\newcommand{\eex}{ \end{ex}}

\newcommand{\bdf}{\begin{definition}}
\newcommand{\edf}{ \end{definition}}

\newcommand{\brem}{\begin{rem}}
\newcommand{\erem}{ \end{rem}}

\newcommand{\blem}{\begin{lemma}}
\newcommand{\elem}{\end{lemma}}

\newcommand{\bhis}{\begin{his}}
\newcommand{\ehis}{ \end{his}}

\newcommand{\bobs}{\begin{observation}}
\newcommand{\eobs}{ \end{observation}}

\newcommand{\bproof}{\begin{proof}}
\newcommand{\eproof}{\end{proof}}

\newcommand{\bexr}{\begin{exer}}
\newcommand{\eexr}{ \end{exer}}

\newcommand{\bnotation}{\begin{notation}}
\newcommand{\enotation}{\end{notation}}

\newcommand{\bconv}{\begin{conv}}
\newcommand{\econv}{ \end{conv}}

\newcommand{\bapp}{\begin{app}}
\newcommand{\eapp}{ \end{app}}

\newcommand{\bfact}{\begin{Observation}}
\newcommand{\efact}{ \end{Observation}}

\newcommand{\bcons}{\begin{cons}}
\newcommand{\econs}{ \end{cons}}

\newcommand{\bdiss}{\begin{diss}}
\newcommand{\ediss}{ \end{diss}}

\newcommand{\Mod}{\mathsf{Mod}}
\renewcommand{\mod}{\mathsf{mod}}

\newcommand{\Add}{\mathsf{Add}}
\newcommand{\add}{\mathsf{add}}

\newcommand{\Proj}{\mathcal{P}}
\newcommand{\proj}{\mathcal{P}^{\mathrm{fin}}}
\newcommand{\GP}{\mathcal{GP}}
\newcommand{\Gp}{\mathcal{GP}^\mathrm{fin}}

\newcommand{\fin}{\mathrm{fin}}

\newcommand{\Filt}{\mathsf{Filt}}
\newcommand{\filt}{\mathsf{filt}}

\newcommand{\arrow}{\longrightarrow}

\renewcommand{\hom}{\operatorname{Hom}}

\newcommand{\Ext}{\operatorname{Ext}}

\renewcommand{\ker}{\operatorname{Ker}}
\newcommand{\coker}{\operatorname{Coker}}
\newcommand{\im}{\operatorname{Im}}

\newcommand{\FPD}{\operatorname{FPD}}
\newcommand{\fpd}{\operatorname{fpd}}
\newcommand{\pd}{\operatorname{pd}}

\newcommand{\Gpd}{\operatorname{Gpd}}

\newcommand{\hookarrow}{\xymatrix{  \ar@{^{(}->}[r] & }}

\title{{\huge Finitistic Dimension Conjectures via Gorenstein Projective Dimension}}

\author{Pooyan Moradifar\footnote{pmoradifar@ut.ac.ir} \\
{\footnotesize \textsc{School of Mathematics, Statistics and Computer Science}} \\ 
{\footnotesize \textsc{ University of Tehran}} \\ 
{\footnotesize \textsc{Tehran, Iran}} 
\and 
Jan \v{S}aroch\footnote{saroch@karlin.mff.cuni.cz} \\
{\footnotesize \textsc{Faculty of Mathematics and Physics, Department of Algebra}} \\ 
{\footnotesize \textsc{Charles University}} \\ {\small \textsc{Sokolovsk\'{a} 83}} \\ 
{\footnotesize \textsc{186 75 Praha 8, Czech Republic}}}
\date{}
\begin{document}
\maketitle
\begin{abstract}
It is a well-known result of Auslander and Reiten that contravariant finiteness of the class $\mathcal{P}^{\mathrm{fin}}_\infty$ (of finitely generated modules of finite projective dimension) over an Artin algebra is a sufficient condition for validity of finitistic dimension conjectures. Motivated by the fact that finitistic dimensions of an algebra can alternatively be computed by Gorenstein projective dimension, in this work we examine the Gorenstein counterpart of Auslander--Reiten condition, namely contravariant finiteness of the class $\mathcal{GP}^{\mathrm{fin}}_\infty$ (of finitely generated modules of finite Gorenstein projective dimension), and its relation to validity of finitistic dimension conjectures. It is proved that contravariant finiteness of the class $\mathcal{GP}^{\mathrm{fin}}_\infty$ implies validity of the second finitistic dimension conjecture over left artinian rings. In the more special setting of Artin algebras, however, it is proved that the Auslander--Reiten sufficient condition and its Gorenstein counterpart are virtually equivalent in the sense that contravariant finiteness of the class $\mathcal{GP}^{\mathrm{fin}}_\infty$ implies contravariant finiteness of the class $\mathcal{P}^{\mathrm{fin}}_\infty$ over any Artin algebra, and the converse holds for Artin algebras over which the class $\mathcal{GP}^{\mathrm{fin}}_0$ (of finitely generated Gorenstein projective modules) is contravariantly finite.\\

\noindent {\footnotesize \textbf{Keywords:}  Contravariant finiteness; Cotorsion pair; Finitistic dimensions; Gorenstein projective dimension; Tilting module.}
\end{abstract}
\pagebreak
\section*{Introduction}
\addcontentsline{toc}{section}{Introduction}

A key problem in homological theory of rings and modules is to understand the range of homological dimensions associated to rings and their modules.
In particular, many problems and (homological) conjectures about rings and modules relate in one way or another to the question of 
\textsl{how large can the projective dimension of modules over a ring be?}
To address this question two classes of modules, namely
the class $\Proj_\infty$ of modules of finite projective dimension and the class $\proj_\infty$ of finitely generated modules of finite projective dimension,
come naturally to the fore.
In order to measure the range of projective dimension of modules in the two classes, the so-called \emph{finitistic dimensions} are natural invariants to consider.
Recall that for any ring $\Lambda$, the \emph{big finitistic dimension} of $\Lambda$ is defined as
\[ \FPD (\Lambda) := \sup \big \{ \pd_\Lambda (M) \mid M \in \Proj_\infty \big \} \: , \]
and the \emph{little finitistic dimension} of $\Lambda$ is defined as
\[ \fpd (\Lambda) := \sup \big \{ \pd_\Lambda (M) \mid M \in \proj_\infty \big \} \: . \] 
The importance of these invariants was emphasized in the 1960s by Bass~\cite{bass.finitistic},
where the following problems due to Rosenberg and Zelinsky, later known as \emph{finitistic dimension conjectures}, 
were advertised~\cite[page~487]{bass.finitistic}:
\begin{description}
\item[First Finitistic Dimension Conjecture:] $\FPD (\Lambda) = \fpd (\Lambda)$.

\item[Second Finitistic Dimension Conjecture:] $\fpd (\Lambda) < + \infty$.
\end{description}
Over the years, the conjectures have been studied mainly in two different but closely related areas, 
namely commutative algebra and representation theory of Artin algebras. In commutative algebra, the conjectures are quite well-understood and they both \emph{fail} 
in general: it is known through the works of Auslander and Buchsbaum~\cite{AB.local}, Bass~\cite{bass.injdim}, and Gruson and Raynaud~\cite{Gruson.Raynaud} that 
if $\Lambda$ is a commutative noetherian local ring, then $\fpd (\Lambda) < + \infty$ and the equality $\FPD (\Lambda) = \fpd (\Lambda)$ holds if 
and only if $\Lambda$ is Cohen-Macaulay. Furthermore, there are examples of commutative noetherian rings $\Lambda$ with $\fpd (\Lambda)=+\infty$; 
see e.g.~\cite[page~276]{Rosen}.
In representation theory of Artin algebras, the finitistic dimension conjectures are not as much well-understood as in the commutative setting:
It is proved by Huisgen-Zimmermann~\cite{Z.domino} that the first finitistic dimension conjecture \textit{fails} in general for a monomial-relation algebra and
examples due to Smal{\o}~\cite{smalo.fin} show that the difference between the first and the second finitistic dimension can indeed be arbitrarily large. 
However, it is still of particular interest to know for which classes of algebras the equality $\FPD (\Lambda) = \fpd (\Lambda)$ holds.
The second  finitistic dimension conjecture is still open in general, but the conjecture is verified for many classes of algebras including algebras of
finite representation type, monomial-relation algebras, radical square/cube zero algebras; see~\cite{Z.tale} for more information in this regard.

It is well-known that understanding the structure of modules in the classes $\proj_\infty$ and $\Proj_\infty$ provides insight to finitistic dimension conjectures,
and for this purpose ``approximation theory'' and ``tilting theory'' turn out to be invaluable tools; see e.g.~\ref{thm:tilting.fpd}. 
The first result via this approach was obtained by Auslander and Reiten~\cite{AR.applications}, who proved that ``contravariant finiteness of the
class $\proj_\infty$'', referred to as the \emph{Auslander--Reiten condition} from now on, is a sufficient condition for validity of the second finitistic dimension conjecture over an Artin algebra. This result was further strengthened or generalized later:
\begin{enumerate} \renewcommand{\theenumi}{(\Roman{enumi})}
\item It was proved by Huisgen-Zimmermann and Smal{\o}~\cite{zimm.smalo} that contravariant finiteness of the class $\proj_\infty$ actually implies
that any module in $\Proj_\infty$ is a direct limit of modules in $\proj_\infty$, and  as a result both finitistic dimension conjectures hold in this case, i.e.\ $\FPD (\Lambda) = \fpd (\Lambda) < + \infty$. Thus, Auslander--Reiten condition is actually a sufficient condition for validity of both finitistic dimension conjectures.
The condition is not, however, necessary as examples due to Igusa, Sma{\o} and Todorov~\cite{Igusa.Todorov} show.

\item  Trlifaj~\cite{trlifaj.approx.little} proved, using tools of approximation theory of modules, that
contravariant finiteness of the class $\proj_\infty$ is still sufficient
for validity of the second finitistic dimension conjecture over left artinian rings.

\item In~\cite{trlifaj.fpd}, Angeleri-H\"{u}gel and Trlifaj presented a somewhat more conceptual proof of the implication
\[ \text{$\proj_\infty$ is contravariantly finite} \implies \FPD (\Lambda) = \fpd (\Lambda) < + \infty \]
mentioned in~(I) using tilting theory, by showing that the class $\proj_\infty$ is contravariantly finite over an Artin algebra if and only
if the cotorsion pair generated by $\proj_\infty$ is induced by a~finitely generated tilting module which renders the
equality $\FPD (\Lambda) = \fpd (\Lambda) < +\infty$; see~\ref{thm:tilting.fpd}.
\end{enumerate}
The above-mentioned results can be recapitulated in the following diagram:
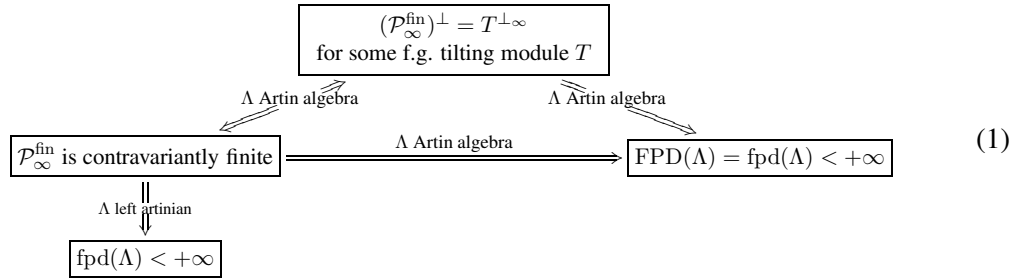
\begin{equation} \label{eq:AR}
\begin{array}{c}
\resizebox{0.8\textwidth}{!}{$
\xymatrixcolsep{0.1cm}
\xymatrix{
& \boxed{{\begin{array}{c}
\text{$(\proj_\infty)^\perp = T^{\perp_\infty}$} \\
\text{for some f.g. tilting module $T$}
\end{array}}} \ar@{=>}[dr]|{\text{{\footnotesize $\Lambda$ Artin algebra}}}  & \\
\boxed{\text{$\proj_\infty$ is contravariantly finite}} \ar@{<=>}[ur]|{\text{{\footnotesize $\Lambda$ Artin algebra}}}
\ar@{=>}[rr]^-{\text{{\footnotesize $\Lambda$ Artin algebra}}} \ar@{=>}[d]|{\text{$\Lambda$ left artinian}} &&
\boxed{\FPD (\Lambda) = \fpd (\Lambda) < +\infty} \\
\boxed{\fpd (\Lambda) < + \infty}  & &
}$}
\end{array}
\end{equation}
The point of departure in the present work is that in studying finitistic dimensions, it is sometimes more convenient to look at some alternative classes of modules,
other than the obvious classes $\proj_\infty$ and $\Proj_\infty$, and two such alternative classes are:
\begin{itemize}
\item the class $\GP_\infty$ of modules of finite Gorenstein projective dimension, and

\item the class $\Gp_\infty$ of finitely generated modules of finite Gorenstein projective dimension.
\end{itemize}
These classes of modules are usually regarded as ``Gorenstein counterparts'' of $\Proj_\infty$ and $\proj_\infty$ in the so-called ``Gorenstein homological algebra'', a branch of homological algebra where the focus is on studying ``Gorenstein modules'' and their respective dimensions; cf.~\cite{enochs.book},~\cite{Gdimensions},~\cite{wang.foundations} and~\cite{iacob.Gorenstein}.
By a well-known result of Holm~\cite[Theorem~2.28]{holm.gorenstein}, finitistic dimensions of a ring can be computed by
modules in $\GP_\infty$ and $\Gp_\infty$. 
More precisely, for \textit{any} ring $\Lambda$,
\begin{equation} \label{holm.FPD}
\FPD (\Lambda)  = \sup \big \{ \Gpd_\Lambda (M) \mid M \in \GP_\infty \big \} 
\end{equation}
and if $\Lambda$ is left notherian, then we also have
\begin{equation} \label{holm.fpd}
\fpd (\Lambda)  = \sup \big \{ \Gpd_\Lambda (M) \mid M \in \Gp_\infty \big \} \: .
\end{equation}
To illustrate the point that using $\GP_\infty$ and $\Gp_\infty$ in place of the obvious classes $\Proj_\infty$ and $\proj_\infty$ is sometimes more convenient,
note that validity of finitistic dimension conjectures for Iwanaga-Gorenstein rings follows almost immediately from~\eqref{holm.FPD} and~\eqref{holm.fpd},
because Iwanaga-Gorenstein rings can be characterized in terms of global finiteness of Gorenstein projective dimension of modules~\cite[Corollary~12.3.2]{enochs.book}; also cf.~\cite{trlifaj.Gorenstein} and~\cite[Corollary~5.2]{bazzoni.tilting}.

\paragraph{Main Problem and Summary of Results.}
In view of the considerations above, it is natural to ask how contravariant finiteness of $\Gp_\infty$, regarded as
the ``Gorenstein counterpart'' of the Auslander--Reiten condition, fits in with the implications mentioned in~\eqref{eq:AR}.
More precisely, our goal in this paper is to investigate the relation between:
\begin{itemize}
\item[(a)] contravariant finiteness of $\Gp_\infty$,

\item[(b)] contravariant finiteness of $\proj_\infty$, and

\item[(c)] validity of finitistic dimension conjectures.
\end{itemize}
It is easy to see in the first place that the condition~(a) above implies validity of the second finitistic dimension conjecture over arin algebras, 
using Eq.~\eqref{holm.fpd} in conjunction with a well-known result of Auslander and Reiten~\cite[Proposition~3.8]{AR.applications}; see also~\cite[Proposition 4.8]{xiIII}. 
Nevertheless, we provide a more complete picture by  showing that ``contravariant finiteness of $\Gp_\infty$'' fits in with the implications in diagram~\eqref{eq:AR} as follows:
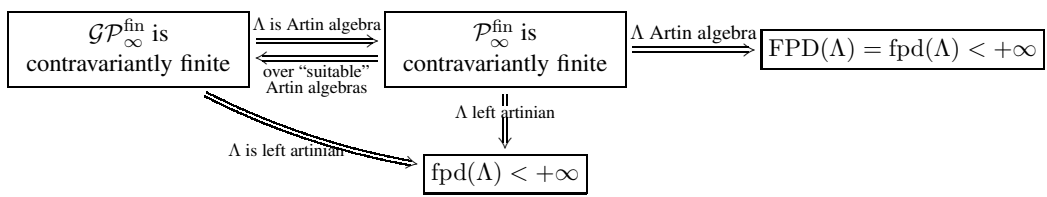
\[
\begin{array}{c}
\resizebox{\textwidth}{!}{$
\xymatrix{
{\boxed{\begin{array}{c}
\text{$\Gp_\infty$ is} \\
\text{contravariantly finite}
\end{array}}} 
\ar@<0.8ex>@{=>}[rr]^-{\text{$\Lambda$ is Artin algebra}} \ar@/_1pc/@{=>}[drr]_-{\text{$\Lambda$ is left artinian}} & &
{\boxed{\begin{array}{c}
\text{$\proj_\infty$ is} \\
\text{contravariantly finite}
\end{array}}}
\ar@<0.8ex>@{=>}[ll]^-{{ \begin{smallmatrix}
\text{over ``suitable''} \\
\text{Artin algebras}
\end{smallmatrix} }}
\ar@{=>}[rr]^-{\text{{\footnotesize $\Lambda$ Artin algebra}}} \ar@{=>}[d]|{\text{$\Lambda$ left artinian}} & &
\boxed{\FPD (\Lambda) = \fpd (\Lambda) < +\infty} \\
& & \boxed{\fpd (\Lambda) < + \infty}  & & &
}$}
\end{array}
\]
The above implications are proved in Section~2 after the preliminary Section~1.
The implication
\[\text{$\Gp_\infty$ is contravariantly finite} \implies \text{$\proj_\infty$ is contravariantly finite}\]
is proved in Theorem~\ref{Main.thm.I} using tilting theory. 
Indeed, it is shown using Proposition~\ref{prop:tilting-like} that when $\Gp_\infty$ is contravariantly finite,
the cotorsion pair generated by $\Gp_\infty$ has a ``tilting-like structure'' in the sense of Definition~\ref{df:tilting.like}, 
and the underlying tilting module can be taken finitely generated which renders $\proj_\infty$ to be contravariantly finite.
The reverse implication, namely
\[ \text{$\proj_\infty$ is contravariantly finite} \implies \text{$\Gp_\infty$ is contravariantly finite} \: , \]
it proved in Theorem~\ref{Main.thmII} for Artin algebras over which the class $\Gp_0$ is contravariantly finite.
Typical examples of such Artin algebras are CM-finite and virtually Gorenstein Artin algebras; cf. Remark~\ref{rem:vGor}.
Thus, one can say that Auslander--Reiten condition and its Gorenstein counterpart are ``almost equivalent''; cf.~\cite[Proposition 4.8]{xiIII}.
In the end, we shift our focus from Artin algebras to the slightly more general setting of left artinian rings and we prove in Theorem~\ref{Main.thmIII} 
that contravariant finiteness of $\Gp_\infty$ still implies validity of the second finitistic dimension conjecture for left artinian rings.
This can be regarded as the ``Gorenstein counterpart'' of the main result of~\cite{trlifaj.approx.little}.

\section{Preliminaries}

\bdes{General Notations, Notions and Conventions}
Throughout the paper, by a ``ring'' we mean an arbitrary non-trivial unital ring.
Such a ring is denoted by $\Lambda$.
We often assume that $\Lambda$ is an \emph{Artin algebra} which means that $\Lambda$ is an algebra over a commutative artinian ring $R$ and
$\Lambda$ is finitely generated as a module over $R$.
By a ``module over $\Lambda$'' or a ``$\Lambda$-module'' we always mean a \textit{left} $\Lambda$-module. 
If $\xymatrixcolsep{0.5cm}\xymatrix{P_\bullet = \cdots \ar[r]^-{f_2} & P_1 \ar[r]^-{f_1} & P_0 \ar[r]^-{f_0} & M \ar[r] & o}$ is
a projective resolution of a $\Lambda$-module $M$, then for any $i \geq 0$ the module $\im (f_i)$ is called the \emph{$i$-th syzygy module} of $M$ in the projective resolution $P_\bullet$.
The class of all $\Lambda$-modules is denoted by $\Mod (\Lambda)$ and the class of \textit{strongly finitely presented modules}, i.e.\ modules with a degreewise finitely generated projective resolution, is denoted by $\mod (\Lambda)$. 
Furthermore, given a class $\mathcal{C}$ of $\Lambda$-modules we let $\mathcal{C}^\mathrm{fin} : = \mathcal{C} \cap \mod (\Lambda)$.

For any integer $n \geq 0$, the class of $\Lambda$-modules of projective dimension at most $n$ is denoted by $\mathcal{P}_n$, and
the class of all $\Lambda$-modules of finite projective dimension is denoted by $\mathcal{P}_\infty$.
These are typical examples of \emph{resolving classes} of modules, that is by definition 
extension closed classes of modules containing $\Proj_0$ which are closed under kernels of epimorphisms.
Dually, a class $\mathcal{C}$ of $\Lambda$-modules is called \emph{coresolving}  if it contains  all injective modules, and it is closed under extensions and 
cokernels of monomorphisms.
Resolving classes are always \emph{syzygy-closed} in the sense that they contain syzygies of projective resolutions of 
their elements. The notions of a ``(co)resolving class'' and a ``syzygy-closed class'' can be defined within $\mod (\Lambda)$ with the obvious modifications. 

$\Ext$-perpendicular classes will be of frequent use in the sequel. 
For any class $\mathcal{C}$ of $\Lambda$-modules let
\begin{align*}
\mathcal{C}^\perp & := \big \{ M \in \Mod (\Lambda) \mid \text{$\Ext^1_\Lambda (C , M) = o$ for all $C \in \mathcal{C}$}  \big \} \:  , \\
{}^\perp \mathcal{C} & := \big \{ M \in \Mod (\Lambda) \mid \text{$\Ext^1_\Lambda (M , C) = o$ for all $C \in \mathcal{C}$} \big \} \: , \\
\mathcal{C}^{\perp_\infty} & := \big \{ M \in \Mod (\Lambda) \mid \text{$\Ext^{\geqslant 1}_\Lambda (C , M) = o$ for all $C \in \mathcal{C}$} \big \} \:  , \\
{}^{\perp_\infty} \mathcal{C} & := \big \{ M \in \Mod (\Lambda) \mid \text{$\Ext^{\geqslant 1}_\Lambda (M , C) = o$ for all $C \in \mathcal{C}$} \big \} \: .
\end{align*}
For any module $M$ we let $M^\perp := \{M\}^\perp$ for simplicity, and  a similar notation is adopted for other  $\Ext$-perpendicular classes of a singleton.
It is easy to see that the class ${}^{\perp_\infty} \mathcal{C}$ is always resolving and the equality $\mathcal{C}^\perp = \mathcal{C}^{\perp_\infty}$ 
holds provided that $\mathcal{C}$ is syzygy-closed.
\edes

\bdes{Filtrations} \label{df:filtrations}
A family $\{ M_\alpha \}_{\alpha \leq \sigma}$ (indexed by an ordinal $\sigma$) is called  a \emph{continuous chain} 
if the inclusion $M_\alpha \subseteq M_{\alpha+1}$ holds for any ordinal $\alpha < \sigma$, and the equality 
$M_\alpha = \bigcup_{\beta < \alpha} M_\beta$ holds for any limit ordinal $\alpha \leq \sigma$.

Given a  class $\mathcal{C}$ of $\Lambda$-modules, a $\Lambda$-module $M$ is called \emph{$\mathcal{C}$-filtered} if there is a continuous 
chain $\{M_\alpha\}_{\alpha \leq \sigma}$ of $\Lambda$-modules such that $M_0 = o$, $M_\sigma = M$, 
and the successive factors $M_{\alpha + 1} / M_{\alpha}$ are isomorphic to an element in $\mathcal{C}$ for any $\alpha < \sigma$.
In this case the family $\{M_\alpha\}_{\alpha \leq \sigma}$ is called a \emph{$\mathcal{C}$-filtration} of $M$ of the length $\sigma$, 
and a $\mathcal{C}$-filtration is said to be \emph{finite} if $\sigma$ is a finite ordinal, i.e.\ a natural number. 
In this case, the module $M$ is said to be \emph{finitely $\mathcal{C}$-filtered}.
\edes

\bnotation
Given a class $\mathcal{C}$ of $\Lambda$-modules, the class of all $\mathcal{C}$-filtered modules is denoted by $\Filt (\mathcal{C})$, 
and the class of of all finitely $\mathcal{C}$-filtered modules is denoted by $\filt (\mathcal{C})$.
\enotation

The following lemma about the length of filtrations of finitely generated modules will be useful later.

\blem \label{ffilt}
Let $\mathcal{C}$ be a class consisting of finitely presented $\Lambda$-modules. 
If $\{M_\alpha\}_{\alpha \leq \sigma}$ 
is a \emph{strict} $\mathcal{C}$-filtration (i.e. $M_\alpha \subsetneqq M_{\alpha+1}$ for every $\alpha < \sigma$) of a finitely generated $\Lambda$-module $M$, then $\sigma$ is finite.
\elem

\bproof
For the sake of contradiction assume that $\sigma$ is infinite. Then we can write
$\sigma = \tau + n$ where $\tau$ is a limit ordinal and $n \geq 0$ is a natural number.
Since $M_\sigma = M$ is finitely generated and $\mathcal{C}$ consists of finitely presented modules, one can deduce inductively that $M_\tau = \bigcup_{\alpha < \tau} M_\alpha$ is finitely generated. 
Thus there exists $\alpha < \tau$ such that $M_\alpha$ contains a finite generating set of $M_\tau$, i.e.\ $M_\tau \subseteq M_\alpha$ which is not possible as the
$\mathcal{C}$-filtration is strict.
\eproof

As it was already mentioned in the introduction, the pivotal class of modules in this paper are modules of ``finite Gorenstein projective dimension''
which are modules finitely resolved by ``Gorenstein projective modules''. 

\bdes{Gorenstein Projective Modules} \label{df:GP.module}
A $\Lambda$-module $M$ is called \emph{Gorenstein projective} if it admits a \emph{complete projective resolution}, that is
an exact sequence 
\[\xymatrix{ P_\bullet = \cdots \ar[r] & P_1 \ar[r]^-{f_1} & P_0 \ar[r]^-{f_{0}} & P_{-1} \ar[r] & \cdots }\]
consisting of projective modules such that $M \cong \ker (f_0)$ and that the sequence remains
exact under $\hom_\Lambda (-,P)$ for every $P\in \Proj_0$.
The class of Gorenstein projective $\Lambda$-modules is denoted by $\GP$.
Note that by the symmetry in the definition, all the modules $\ker (f_i)$ in the complete projective resolution $P_\bullet$ are also Gorenstein projective.
In the special case where all $f_i$ and $P_i$ are equal in $P_\bullet$, the Gorenstein projective module $M$ is called \emph{strongly Gorenstein projective}.
\edes

The importance of strongly Gorenstein projective modules lies in the following construction from~\cite{bennis.strong}; see also~\cite[Theorem 11.1.12]{wang.foundations}.

\bdes{Construction} \label{construct:S}
Let $M$ be a Gorenstein projective $\Lambda$-module and
\[\xymatrix{
P_\bullet = \cdots \ar[r]^-{d_2} & P_1 \ar[r]^-{d_1} & P_0 \ar[r]^-{d_0} & P_{-1} \ar[r]^-{d_{-1}} & \cdots
}\]
be a complete projective resolution of $M$.
Let $P : = \bigoplus_{n \in \mathbb{Z}} P_n$ and $\partial : P \arrow P$ be the $\Lambda$-homomorphism induced by $d_i$, i.e.\ $\partial \restriction P_i = d_i$ for every $i \in \mathbb{Z}$.
It is then easy to see (see~\cite[Theorem 2.7]{bennis.strong} or \cite[Theorem 11.1.12]{wang.foundations}) that $S := \ker (\partial)$ is a strongly Gorenstein projective module containing $M$ as a direct summand.
\edes

\bdes{Gorenstein Projective Dimension} \label{df:G-dimensions}
It is well-known that the class $\GP$ is resolving~\cite[Theorem~2.5]{holm.gorenstein} and hence one can define \emph{Gorenstein projective dimension} of modules by
resolving modules by the class $\GP$; cf.~\cite{holm.gorenstein} and~\cite{enochs.book} for more information.  
The Gorenstein projective dimension of $\Lambda$-modules is denoted by $\Gpd_\Lambda (-)$.
For any integer $n \geq 0$ the class of $\Lambda$-modules of Gorenstein projective dimension at most $n$ is denoted by $\GP_n$.
We also let $\GP_\infty$ be the class of  $\Lambda$-modules of finite Gorenstein projective dimension.
\edes

\brem \label{rem:GP}
It is well-known that $\GP_n$ has two-of-three property~\cite[Section~11.3]{wang.foundations} and it is closed under filtrations~\cite[Theorem~3.4]{enochs.transext}. 
\erem

Modules of finite projective dimension are  important partly because they behave more or less similarly to modules over rings of finite global dimension.
On the other hand, the ``Gorenstein version'' of these modules, namely modules of finite Gorenstein projective dimension,
have historically been studied as the class of modules which behave more or less similarly to 
modules over Gorenstein rings; cf.~\cite{beyond.and.back},~\cite{Gdimensions} and~\cite{iacob.Gorenstein}. 

The following classical result~\cite[Theorem 11.3]{wang.foundations} shows that Gorenstein projective dimension of modules in $\GP_\infty$
can be measured via vanishing of $\Ext$-functors.

\bthm \label{thm:GPT}
Let $\Lambda$ be a ring and $n \geq 0$ be an integer. 
The following statements are equivalent for any $M \in \GP_\infty$:
\begin{enumerate}
\item $\Gpd_\Lambda (M) \leq n$;

\item $\Ext^{n+i}_\Lambda (M , N)$ for all $i \geq 1$ and $N \in \Proj_\infty$.

\item $\Ext^{n+1}_\Lambda (M , N)$ for all $N \in \Proj_\infty$.
\end{enumerate}
\ethm

The following lemma will be used later in the proofs of a couple of results.

\blem \label{lem:LMN}
For any integer $n \geq 0$ and $M \in \Gp_n$, there exists a short exact sequence
\[\xymatrix{
o \ar[r] & M \ar[r] & P \ar[r] & G \ar[r] & o
}\]
where $P \in \mathcal{P}_n^\mathrm{fin}$ and $G \in \Gp_0$.
\elem

\bproof
We proceed by induction on $n$.
For $n=0$ there exists---essentially by definition of a Gorenstein projective module---a short exact sequence
\[\xymatrix{o \ar[r] & M \ar[r]^-u & F \ar[r]^-v & C \ar[r] & o }\]
of $\Lambda$-modules where $F$ is projective and $C$ is Gorenstein projective. 
Since $F$ is a direct summand of a free module, we may add a suitable projective summand to 
$\xymatrix{F \ar[r]^-v & C}$ so that we can assume $F$ is free. 
Now since $M$ is finitely generated,
$u : M \arrow F$ factors through some finitely generated free $\Lambda$-module $\Lambda^n$ and thence we obtain a short exact sequence
\[\xymatrix{o \ar[r] & M \ar[r]^-f & \Lambda^n \ar[r]^-g & N \ar[r] & o }\]
wherein $N \in \GP_\infty$ by Remark~\ref{rem:GP}, and $\hom_\Lambda(f,P)$ is surjective for each projective module $P$ because $\hom_\Lambda(u,P)$ was such.
Consequently $\Ext^1_\Lambda(N,P) = 0$ for every $P \in \Proj_0$ and then Theorem~\ref{thm:GPT} yields $N \in \GP_0$.
Assume now that $n \geq 1$ and the assertion holds for all modules in $\Gp_{n-1}$.
Consider a short exact sequence
\[\xymatrix{o \ar[r] & M \ar[r] & F \ar[r] & N \ar[r] & o}\]
where $F$ is free, and note that $M \in \Gp_{n-1}$.
Thus, by the inductive hypothesis, we have a short exact sequence
\[\xymatrix{
o \ar[r] & M \ar[r] & Q \ar[r] & H \ar[r] & o
}\]
with $Q \in \mathcal{P}_{n-1}^\mathrm{fin}$ and $H \in \Gp_0$. 
Forming the pushout of the last two short exact sequences, we obtain a short exact sequence
\[\xymatrix{
o \ar[r] & Q \ar[r] & U \ar[r] & N \ar[r] & o
}\]
where $U\in\Gp_0$.
Finally, we form one more pushout
\[\xymatrix{
&  & o \ar[d] & o \ar@{..>}[d] &  \\
o \ar[r] & Q \ar[r] \ar@{=}[d] & U \ar[r] \ar[d] & N \ar@{..>}[d] \ar[r] & o \\
o \ar@{..>}[r] & Q \ar@{..>}[r] & \Lambda^n \ar@{..>}[r] \ar[d] & P \ar@{..>}[r] \ar@{..>}[d] & o  \\
& & G \ar@{=}[r] \ar[d] & G \ar@{..>}[d] & \\
& & o & o & 
}\]
where we use the already proved ``case $n=0$'' to obtain the middle column with $G\in \Gp_0$. 
The short exact sequence on the right-hand side column is obviously the desired one.
\eproof

Next, we review some definitions and facts from approximation theory of modules; we refer to~\cite{GT} for more information and proofs of the standard facts mentioned below.

\bdes{Approximations}
Given a $\Lambda$-module $M$ and a class $\mathcal{C}$ of $\Lambda$-modules, a $\Lambda$-homomorphism $f : C \arrow M$ with $C \in \mathcal{C}$ is said to be a \emph{$\mathcal{C}$-precover (of $M$)}
if any $\Lambda$-homomorphism $g : C' \arrow M$ with $C' \in \mathcal{C}$ factors through $f$.
\[\xymatrix{
& C' \ar[d]^-g \ar@{-->}[dl] \\
C \ar[r]_-f & M
}\]
The class $\mathcal{C}$ is called \emph{precovering} if any $\Lambda$-module has a $\mathcal{C}$-precover.
The dual of the notion ``precover'' is called a \emph{preenvelope} and subsequently we may speak of 
a \emph{preenveloping class} of $\Lambda$-modules; cf.~\cite{GT}.

Approximation theory of modules, also known as the ``theory of covers and envelopes'', originates in the work of Enochs~\cite{enochs.coversenv}
on torsion-free and flat covers of modules, and also earlier work of Auslander and Smal{\o}~\cite{aus.preproj,ARseq.in.subcat} in the realm of
 representation theory of Artin algebras.
In the latter setting, the related notions are often confined to the class of finitely generated modules over Artin algebras, and
in this setting the term \emph{contravariantly finite} is often  synonymously used instead of ``precovering in $\mod (\Lambda)$''.

A typical situation where a $\Lambda$-homomorphism $f : C \arrow M$ with $C \in \mathcal{C}$ happens to be a $\mathcal{C}$-precover is 
when $f$ is surjective and $\ker (f) \in \mathcal{C}^\perp$. In this case, $f : C \arrow M$ is called a \emph{special $\mathcal{C}$-precover (of $M$)}.
Dually, a $\Lambda$-homomorphism $f : M \arrow C$ with $C \in \mathcal{C}$ is called a \emph{special $\mathcal{C}$-preenvelope (of $M$)} if $f$ is injective and $\coker (f) \in {}^\perp \mathcal{C}$. Accordingly, the class $\mathcal{C}$ is called \emph{special precovering} (respectively, \emph{special preenveloping}) if any $\Lambda$-module has a special $\mathcal{C}$-precover (respectively, special $\mathcal{C}$-precover).
\edes

As the following proposition shows, in the setting of finitely generated modules over Artin algebras,  contravariantly finite classes with suitable closure properties actually provide for special precovers.

\bprop \label{prop:approx.in.mod}
Let $\Lambda$ be an Artin algebra and $\mathcal{C}$ be a
resolving class in $\mod (\Lambda)$. If $\mathcal{C}$ is contravariantly finite, then $\mathcal{C}$ is special precovering in $\mod (\Lambda)$.
\eprop

\bproof
By the hypothesis each finitely generated $\Lambda$-module $M$ has a $\mathcal{C}$-precover $f : C \arrow M$ in $\mod (\Lambda)$ which is surjective because $\Lambda \in \mathcal{C}$.
Furthermore, $f$ can be taken to be ``left minimal'' 
in the sense that $f$ cannot factor through a $\Lambda$-homomorphism $\alpha : C \arrow C$ unless $\alpha$ is an isomorphism; cf.~\cite[Corollary 5.10]{GT} or~\cite[Theorem 2.4 ]{book}.
Now Wakamatsu Lemma~\cite[Lemma 5.13]{GT} in $\mod (\Lambda)$---see also~\cite[Lemma~1.3]{AR.applications}---implies that such a 
left minimal $\mathcal{C}$-precover in $\mod (\Lambda)$ is a special $\mathcal{C}$-precover in $\mod (\Lambda)$.
\eproof

A useful machinery to systematically detect or construct classes of modules which are special precovering or special preenveloping
is the notion of a ``cotorsion pair''. 

\bdes{Cotorsion Pairs} \label{df:cotorsion.pair}
A  pair $(\mathcal{L} , \mathcal{R})$ of two classes of $\Lambda$-modules is said to be a \emph{cotorsion pair} 
if $\mathcal{L}={}^\perp \mathcal{R}$ and $\mathcal{R} = \mathcal{L}^\perp$. 
Given a class $\mathcal{C}$ of $\Lambda$-modules, it is easily seen that $\big ( {}^\perp (\mathcal{C}^\perp) , \mathcal{C}^\perp \big )$ is a cotorsion pair
called the \emph{cotorsion pair generated by $\mathcal{C}$}.

The components of a cotorsion pair are known to have some dual properties:
A result known as  the Rozas Lemma~\cite[Lemma 5.24]{GT} states that in a cotorsion pair $\mathfrak{C}:=(\mathcal{L} , \mathcal{R})$ the left-hand side class $\mathcal{L}$ is resolving if and only if the right-hand side class $\mathcal{R}$ is coresolving. 
In this case the cotorsion pair $\mathfrak{C}$ is said to be \emph{hereditary}.
It is easy to see that every cotorsion pair generated by a syzygy-closed class of modules is hereditary.
Another duality result in cotorsion pairs is Salce Lemma~\cite[Lemma 5.20]{GT} which states that in the cotorsion pair $\mathfrak{C}$ the left-hand side 
class $\mathcal{L}$ is special precovering if and only if the right-hand side class $\mathcal{R}$ is special preenveloping. 
In this case the cotorsion pair $\mathfrak{C}$ is said to be \emph{complete}.
The following key result due to Eklof and Trlifaj~\cite{trlifaj.Extvanish} shows abundance of complete cotorsion pairs; cf.~\cite[Theorem 6.11]{GT}.
\edes

\bthm \label{thm:ET}
If $\mathcal{S}$ is a set of $\Lambda$-modules, then for any $\Lambda$-module $M$ there exists a short exact sequence
\[\xymatrix{ o \ar[r]  & M \ar[r]^f & N \ar[r] & L \ar[r] & o}\]
where $N \in \mathcal{S}^\perp $ and $L$ is $\mathcal{S}$-filtered.
In particular, $f$ is a special $\mathcal{S}^\perp$-preenvelope and the cotorsion pair generated by $\mathcal{S}$ is complete.
\ethm

Theorem~\ref{thm:ET} provides in particular a relatively concrete description of modules in the double-perpendicular 
class ${}^\perp (\mathcal{S}^\perp)$ when $\mathcal{S}$ is a set. 

\bnotation
For a class $\mathcal{C}$ of modules, $\add_\Lambda (\mathcal{C})$ denotes the class of all direct summands of finite direct sums of modules in $\mathcal{C}$.
\enotation

\bcoro \label{cro:ET}
Let $\mathcal{S}$ be a set of $\Lambda$-modules.
\begin{enumerate}
\item the class ${}^\perp (\mathcal{S}^\perp)$ consists precisely of  direct summands of modules filtered by $\mathcal{S} \cup \{\Lambda\}$. 

\item If $\mathcal{S}$ consists of finitely presented $\Lambda$-modules 
and $\Lambda \in \mathcal{S}$, then ${}^\perp (\mathcal{S}^\perp) \cap \mod (\Lambda) = \add_\Lambda \big ( \filt (\mathcal{S}) \big )$.
\end{enumerate}
\ecoro

\bproof
Part~(i) is a classical result, see e.g.~\cite[Corollary~6.13]{GT}. 
In order to prove part~(ii), notice first that the inclusion
$\add_\Lambda \big ( \filt (\mathcal{S}) \big ) \subseteq {}^\perp (\mathcal{S}^\perp) \cap \mod (\Lambda)$ 
can be proved readily---either by a straightforward argument or using part~(i)---and so
it remains to prove the reverse inclusion.
If $M \in {}^\perp (\mathcal{S}^\perp) \cap \mod (\Lambda)$, then  $M$ is a direct summand of some $\mathcal{S}$-filtered module $N$ by part~(i).
Then it can be proved, say by Hill Lemma~\cite[Theorem 7.10]{GT}, that we can replace $N$ with a finitely presented module, and so we can
take $N$ finitely $\mathcal{S}$-filtered by Lemma~\ref{ffilt}.
Consequently, $M \in \add_\Lambda \big ( \filt (\mathcal{S}) \big )$, and this finishes the proof.
\eproof

The double perpendicular class ${}^\perp (\mathcal{S}^\perp)$ assumes a simpler description over left artinian rings as it turns out that
in this case special precovers of simple modules are enough to determine the structure of all the modules in ${}^\perp (\mathcal{S}^\perp)$.

\bprop \label{prop:artinian.double.perp}
Let $\Lambda$ be a left artinian ring and $\mathcal{S}$ be a set $\Lambda$-modules.
Let $\{S_1 , \ldots , S_n \}$ be a complete set of simple $\Lambda$-modules,
and for every $1 \leq i \leq n$ choose a special ${}^\perp (\mathcal{S}^\perp)$-precover $A_i \arrow S_i$.
Then any (finitely generated) $\Lambda$-module $M$ has special ${}^\perp (S^\perp)$-precover $A \arrow M$ where $A$ is (finitely) filtered by
$\mathcal{C} = \{A_1 , \ldots , A_n \}$.
In particular, the class ${}^\perp (\mathcal{S}^\perp)$ coincides with the class of direct summands of $\mathcal{C}$-filtered modules.
\eprop

\bproof
See~\cite[Corollary 17.19]{GT}.
\eproof

The above results indicate the use of approximation theory in decoding structure of modules. 
Another useful tool serving this purpose is the notion of a ``tilting module''.

\bdes{Tilting Modules} \label{df:tilting}
Let $n \geq 0$ be an integer. A $\Lambda$-module $T$ is  said to be an \emph{$n$-tilting module} if it satisfies the following conditions:
\begin{enumerate}
\item[(T1)] $\pd_\Lambda (T) \leq n$. 
 
\item[(T2)] $\Ext^n_\Lambda (T , T^{(\kappa)})=0$ for any integer $n \geq 1$ and any cardinal $\kappa$. 
Here $T^{(\kappa)}$ denotes the direct sum of $\kappa$ copies of $T$.

\item[(T3)] There exists an exact sequence
\[\xymatrix{ o \ar[r] & \Lambda \ar[r] & T_0 \ar[r] & \cdots \ar[r] & T_m \ar[r] & o}\]
where $T_i \in \Add_\Lambda (T)$ for all $1 \leq i \leq m$.
Here $\Add_\Lambda (T)$ denotes the class of modules isomorphic to a direct summands of $T^{(\kappa)}$ for some cardinal $\kappa$.
\end{enumerate}
In this case, the class $T^{\perp_\infty}$ is called the \emph{tilting class} associated to $T$, and 
the complete and hereditary cotorsion pair $\big ( {}^\perp (T^{\perp_\infty}) , T^{\perp_\infty}\big )$  
is called the \emph{tilting cotorsion pair induced by $T$}.
A~cotorsion pair is said to be an \emph{$n$-tilting cotorsion pair} if it is induced by an $n$-tilting module.
\edes

Tilting modules are the main objects of study in tilting theory with myriads of applications in representation theory of algebras~\cite{tilting.handbook} and structure theory of modules~\cite[Part~III]{GT}.  
Applications of tilting modules to finitistic dimension conjectures were first observed in~\cite{trlifaj.fpd}, where
the authors prove, among other things, the following results; cf.~\cite[Chapter 17]{GT}.

\bthm \label{thm:tilting.fpd}
Let $\mathfrak{P}$ be the cotorsion pair generated by $\proj_\infty$.
\begin{enumerate}
\item If $\Lambda$ is left noetherian, then $\fpd (\Lambda) < +\infty$ if and only if the cotorsion pair $\mathfrak{P}$ is tilting.

\item If $\Lambda$ is an Artin algebra, the class $\proj_\infty$ is contravariantly finite if and only if $\mathfrak{P}$ is a~tilting cotorsion pair induced by a finitely
generated tilting $\Lambda$-module $T$. In this case, $\Proj_\infty = {}^\perp (T^{\perp_\infty})$ and $\FPD (\Lambda) = \fpd (\Lambda)<+\infty$. 
\end{enumerate}
\ethm

Coming  back to our main problem, namely the relation between finitistic dimension conjectures and contravariant finiteness of the class $\Gp_\infty$,
it is natural to consider the cotorsion pair $\mathfrak{G}$ generated by the class $\Gp_\infty$ and employ ideas parallel to~\cite{trlifaj.fpd}, in particular Theorem~\ref{thm:tilting.fpd} mentioned above, to gain insight. 
This approach, however, does not work directly as the cotorsion pair $\mathfrak{G}$ cannot be tilting except when $\Gp_0 = \proj_0$---Artin algebras with this
property are known as \emph{CM-free} in the literature; see e.g.~\cite{chen.gorenstein}.
The remedy to this obstacle is the observation that although the cotorsion pair $\mathfrak{G}$ is almost never tilting, it still has a ``tilting-like structure''
in the sense defined below; see Theorem~\ref{Main.thm.I}.

\bdf \label{df:tilting.like}
A cotorsion pair $(\mathcal{L} , \mathcal{R})$ is said to be \emph{$n$-tilting-like} (for some integer $n \geq 0$) if $\mathcal{R} = (T \oplus S)^{\perp_\infty}$, 
where $T$ is a tilting $\Lambda$-module and $S$ is a strongly Gorenstein projective module.
Needless to say, for $S=o$ we recover the tilting cotorsion pair induced by $T$.
\edf

The following proposition says that in a tilting-like cotorsion pair $(\mathcal{L} , \mathcal{R})$ with the underlying tilting module $T$,
the class ${}^\perp (T^{\perp_\infty})$ coincides with the subclass of $\mathcal{L}$ consisting of modules of finite projective dimension.

\bprop \label{prop:TS} 
Let $ U:= T \oplus S$ where $T$ is an $n$-tilting $\Lambda$-module and $S$ is a strongly Gorenstein projective module. 
If $(\mathcal{L}_U , \mathcal{R}_U)$ is the cotorsion pair with $\mathcal{R}_U = U^{\perp_\infty}$, then
$\mathcal{L}_U \cap \mathcal{P}_n = {}^\perp (T^{\perp_\infty})$.
\eprop

\bproof
The inclusion ${}^\perp (T^{\perp_\infty}) \subseteq \mathcal{L}_U \cap \Proj_\infty$ holds because $T \in \mathcal{L}_U \cap \Proj_n$.
As for the reverse inclusion, let $M \in \mathcal{L}_U \cap \Proj_n$ and $A \in T^{\perp_\infty}$. 
By~\cite[Proposition 13.13]{GT} the module $A$ has an $\Add_\Lambda (T)$-resolution
\[ \xymatrix{
\cdots \ar[r] & T_1 \ar[r]^-{d_1} & T_0 \ar[r]^-{d_0} & A \ar[r] & o } \: . \]
Using condition (T2) in~\ref{df:tilting} and Theorem~\ref{thm:GPT} it is easy to see that the modules $T_i$ in the above sequence belong to $\mathcal{R}_U$ 
and so $\Ext^{\geqslant 1}_\Lambda (M,T_i)=o$ for every $i \geq 0$.
Now since $M \in \Proj_n$, it follows from dimension shifting that
\[ \Ext^1_\Lambda (M,A) \cong \Ext^{n+1}_\Lambda \big (M, \im (d_n) \big ) = o \: . \]
Consequently, $M \in {}^\perp (T^{\perp_\infty})$ and this finishes the proof.
\eproof

Tilting-like cotorsion pairs have recently been studied in~\cite{wang.tilting}, where the authors prove the following
characterization theorem.

\bthm[{\cite[Theorem 1.1]{wang.tilting}}] \label{thm:tilting-like}
Let $\mathfrak{C} := (\mathcal{L} , \mathcal{R} )$ be a hereditary cotorsion pair generated by a set $\mathcal{S}$ of $\Lambda$-modules.
The following statements are equivalent for the cotorsion pair $\mathfrak{C}$:
\begin{enumerate}
\item $\mathfrak{C}$ is $n$-tilting like.

\item There is an $n$-tilting module such that $\mathcal{L} \cap \mathcal{R} = \Add_\Lambda (T)$.

\item $\mathcal{L} \cap \mathcal{R}$ is closed under direct sums and there exists a strongly Gorenstein projective module $S \in \mathcal{L}$ that contains 
some $n$-th syzygy module of every module in $\mathcal{S}$ as a direct summand.

\item $\mathcal{L} \subseteq \GP_n$, $\mathcal{L} \cap \mathcal{R} \subseteq \Proj_n$ and $\mathcal{L} \cap \mathcal{R}$ is closed under direct sums.
\end{enumerate}
\ethm

We will use the following instance of Theorem~\ref{thm:tilting-like} in the next section to prove one of our main theorems, namely Theorem~\ref{Main.thm.I}.

\bprop \label{prop:tilting-like}
Let $\Lambda$ be a ring and $(\mathcal{L} , \mathcal{R})$ be the cotorsion pair generated by a syzygy-closed subclass $\mathcal{C}$ of $\mod (\Lambda)$.
If $\Gp_0 \subseteq \mathcal{C} \subseteq \Gp_n$ for some integer $n \geq 0$, then there exists an $n$-tilting module $T$ and a strongly Gorenstein projective module $S$ such that $\mathcal{R} = (T \oplus S)^{\perp_\infty}$. 
If furthermore $\Lambda$ is an Artin algebra and $\mathcal{L}^\fin$ is contravariantly finite, then the tilting module $T$ can be taken in $\mod (\Lambda)$.
\eprop

\bproof
Replacing $\mathcal{C}$ with a representative set of its element we may assume from the outset that $\mathcal{C}$ is a set.
Since $\mathcal{C}$ consists of strongly finitely presented modules, the class $\mathcal{R} = \mathcal{C}^\perp$ and hence $\mathcal{L} \cap \mathcal{R}$ 
is closed under direst sums.
Furthermore, for every $M \in \mathcal{C}$ it follows from Lemma~\ref{lem:LMN} that any $n$-th syzygy module of 
$M$ has a complete projective resolution whose cycles lie in $\Gp_0 \subseteq \mathcal{L}$.
Thus it follows from Construction~\ref{construct:S} that there exists a strongly Gorenstein projective module $S_M \in \mathcal{L}$
which contains an $n$-th syzygy module of $M$ as a direct summand.
Therefore, it follows from Theorem~\ref{thm:tilting-like} that $\mathcal{R} = (T \oplus S)^{\perp_\infty}$ for some $n$-tilting module $T$ and 
strongly Gorenstein projective module $S$.

As for the second part of the assertion, note that the tilting module $T$
is constructed in general as follows (cf. proof of part~(2)$\implies$(3) in~\cite[Theorem~1.1]{wang.tilting}): 
By considering iterated  special $\mathcal{R}$-preencelopes of $\Lambda$ one can construct an exact sequence
\[\xymatrix{o \ar[r] & \Lambda \ar[r] & T_0 \ar[r] & \cdots \ar[r] & T_n \ar[r] & o}\]
where $T_i \in \mathcal{L} \cap \mathcal{R}$, and it then can be proved that $T:=\bigoplus_{i = 0}^n T_i$ is the desired tilting module.
Now if $\mathcal{L}^\fin$ is contravariantly finite, then it follows from~\cite[Theorem~5.3]{saroch.telescope} that
$(\mathcal{L}^\fin , \mathcal{R}^\fin)$ is a complete cotorsion pair in $\mod (\Lambda)$. So one can repeat the construction process of
$T$ using iterated special $\mathcal{R}^\fin$-preenvelopes of $\Lambda$ inside $\mod (\Lambda)$ and thereby $T$ can be taken
finitely generated.
\eproof

\section{Contravariant finiteness of $\Gp_\infty$}

In this section we prove our main results about the relation between contravariant finiteness of $\Gp_\infty$ and finitistic dimension conjectures, advertised 
earlier in the introduction.
We start with the following theorem which states that the Gorenstein version of Auslander--Reiten condition actually implies the
usual Auslander--Reiten condition; compare~\cite[Proposition 4.8]{xiIII}.

\bthm \label{Main.thm.I}
Let $\Lambda$ be an Artin algebra and consider the following statements about the
cotorsion pair $(\mathcal{L} , \mathcal{R})$ generated by a syzygy-closed class $\mathcal{S} \subseteq \mod (\Lambda)$:
\begin{enumerate}
\item $\mathcal{L}^\fin$ is contravariantly finite;

\item $\mathcal{R} = (T \oplus S)^{\perp_\infty}$ for some finitely generated $n$-tilting $\Lambda$-module $T$ and some strongly Gorenstein projective module $S$.

\item $\mathcal{L}^\fin \cap \mathcal{P}_\infty$ is contravariantly finite.
\end{enumerate}
If $\Gp_0 \subseteq \mathcal{S} \subseteq \Gp_\infty$, then (i)$\implies$(ii)$\implies$(iii).
In particular, contravariant finiteness of $\Gp_\infty$ implies contravariant finiteness of $\proj_\infty$
\ethm

\bproof
Let $\{S_1 , \ldots , S_m\}$ be a complete set of simple $\Lambda$-modules and for every $1 \leq i \leq m$ let $A_i \arrow S_i$ be a 
special ${}^\perp (\mathcal{S}^\perp)$-precover of $S_i$.
Let $n := \sup \{ \Gpd_\Lambda (A_i) \mid 1 \leq i \leq m \}$.
Since the class $\GP_n$ is closed under filtrations~\ref{rem:GP}, it follows from Proposition~\ref{prop:artinian.double.perp}
that $\mathcal{S} \subseteq \mathcal{L} \subseteq \GP_n$.
Now the implication (i)$\implies$(ii) follows from Proposition~\ref{prop:tilting-like}.
As for the implication (ii)$\implies$(iii), it follows from Proposition~\ref{prop:TS} that
$\mathcal{L}^\fin \cap \Proj_\infty = {}^\perp (T^{\perp_\infty}) \cap \mod (\Lambda)$, and it is well-known from classical tilting theory 
that this class is contravariantly finite; cf.~\cite[Lemma 17.26]{GT}.
\eproof

Next we are going to prove in Theorem~\ref{Main.thmII} the converse of the implication
\[\text{$\Gp_\infty$ is contravariantly finite} \implies \text{$\proj_\infty$ is contravariantly finite}\]
for Artin algebras over which $\Gp_0$ is contravariantly finite; cf.~\ref{rem:vGor}.
The key step in the proof is the observation that modules in the class $\Gp_n$ are precisely the modules obtained
as an extension of a module in $\Gp_0$ by a module in $\proj_n$, and that glueing the two classes by extension preserves
contravariant finiteness; see~\ref{prop:star.operation}.
In order to precisely state and prove these considerations, we require some preparatory results.

\bnotation
Let $(\mathcal{U} , \mathcal{V})$ be a pair of subclasses of $\Mod (\Lambda)$. Denote by $\mathcal{U} \star \mathcal{V}$ the class of all modules $M$ which
sit in a short exact sequence of the form
\[\xymatrix{
o \ar[r] & U \ar[r] & M \ar[r] & V \ar[r] & o
}\]
where $U \in \mathcal{U}$ and $V \in \mathcal{V}$.
In other words, $\mathcal{U} \star \mathcal{V}$ is the class of all modules which are an extension of a module in $\mathcal{U}$ by a module in $\mathcal{V}$. 
\enotation

The importance of the operation ``$\star$'' for our purposes lies in the following fact from~\cite{sikko.smalo}; see also~\cite{chen.extensions}.

\bprop \label{prop:star.operation}
Let $\Lambda$ be a ring.
If $\mathcal{U}$ and $\mathcal{V}$ are precovering (resp., preenveloping) classes in $\Mod (\Lambda)$, then then class $\mathcal{U} \star \mathcal{V}$ is also precovering (reps., preenveloping) in $\Mod (\Lambda)$, and this statement remains valid if we replace $\Mod (\Lambda)$ with $\mod (\Lambda)$ in case $\Lambda$ is an Artin algebra.
\eprop

It is clear from the definition that $\mathcal{U} \star \mathcal{V} \subseteq \filt (\mathcal{U} \cup \mathcal{V})$, and the following lemma 
provides a sufficient condition for the equality.

\blem \label{lem:UV.filt}
Let $\Lambda$ be a ring and $(\mathcal{U} , \mathcal{V})$ be a pair of extension closed subclasses of $\Mod (\Lambda)$ which contain the zero 
module. If $\mathcal{U} \subseteq {}^\perp \mathcal{V}$, then $\mathcal{U} \star \mathcal{V} = \filt (\mathcal{U} \cup \mathcal{V})$
\elem

\bproof
We need to show that if $M$ is a finitely $(\mathcal{U} \cup \mathcal{V})$-filtered $\Lambda$-module, 
then there exists a short exact sequence
\[\xymatrix{
o \ar[r] & U \ar[r] & M \ar[r] & V \ar[r] & o
}\]
where $U \in \mathcal{U}$ and $V \in \mathcal{V}$.
By the hypothesis, there exists a finite $(\mathcal{U} \cup \mathcal{V})$-filtration
\[\xymatrix{
o = M_0 \subseteq M_1 \subseteq \cdots \subseteq M_k = M \: .
}\]
We prove the assertion by induction on $k$, the length of the filtration: For $k=1$ the assertion holds trivially. 
Assume that $k>1$ and that the assertion holds for all modules with $(\mathcal{U} \cup \mathcal{V})$-filtration of the length $k-1$. 
Then we have a short exact sequence
\[\xymatrix{
o \ar[r] & L \ar[r]^-f & M \ar[r] & N \ar[r] & o
}\]
where $L$ has a $(\mathcal{U} \cup \mathcal{V})$-filtration of the length $k-1$ and $N \in \mathcal{U} \cup V$. By the induction hypothesis, the module $L$ sits in a short exact sequence
\[\xymatrix{
o \ar[r] & U \ar[r]^-g & L \ar[r] & V \ar[r] & o \: ,
}\]
where $U \in \mathcal{U}$ and $V \in \mathcal{V}$. 
Thus we get the short exact sequence
\begin{equation} \tag{$\ast$}
\xymatrix{
o \ar[r] & U \ar[r]^-{f \circ g} & M \ar[r]^-h & 
W \ar[r] & o \: ,
}\end{equation}
wherein $W$ sits in a short exact sequence of the form
\begin{equation} \tag{$\ast \ast$}
\xymatrix{ o \ar[r] & V \ar[r] & W \ar[r] & N \ar[r] & o } \: .
\end{equation}
Now if $N \in \mathcal{V}$, then $W \in \mathcal{V}$ and ($\ast$) is the desired short exact sequence. Otherwise, $N \in \mathcal{U}$ and so the short exact sequence ($\ast \ast$) splits. In this case, we have also the split short exact sequence
\[\xymatrix{ o \ar[r] & N \ar[r] & W \ar[r]^-\alpha & V \ar[r] & o } \: ,\]
with which we can form the short exact sequence
\begin{equation} \tag{$\dagger$}
\xymatrix{
o \ar[r] & K \ar[r] & M \ar[r]^-{\alpha \circ h} & 
V \ar[r] & o}
\end{equation}
%
wherein $K$ sits in a short exact sequence of the form
\[\xymatrix{
o \ar[r] & U \ar[r] & K \ar[r] & N \ar[r] & o
} \: . \]
and hence it belong to $\mathcal{U}$. 
Now since $\mathcal{U}$ is closed under extensions, we have $K \in \mathcal{U}$ and so the
short exact sequence ($\dagger$) shows that $M \in \mathcal{U} \star \mathcal{V}$.
The proof is thus complete.
\eproof

We are now ready to prove our second main theorem.

\bthm \label{Main.thmII}
If $\Lambda$ is an Artin algebra, then:
\begin{enumerate}
\item For any integer $n \geq 0$ the equality $\Gp_n = \add \big ( \Gp_0 \star \proj_n \big )$ holds.
In particular, if $\proj_n$ and $\Gp_0$ are contravariantly finite, then $\Gp_n$ is contravariantly finite.

\item If $\Gp_0$ and $\proj_\infty$ are contravariantly finite, then $\Gp_\infty$ is contravariantly finite.
\end{enumerate} 
%
\ethm

\bproof
Part~(i): Since the class ${}^\perp \big ( (\Gp_0 \cup \mathcal{P}_n^\mathrm{fin})^\perp \big )$ is resolving, it follows from 
Lemma~\ref{lem:LMN} that $\Gp_n \subseteq {}^\perp \big ( (\Gp_0 \cup \mathcal{P}_n^\mathrm{fin})^\perp \big )$,
and so the equality $\Gp_n = \add \big ( \filt (\Gp_0 \cup \proj_n)  \big )$ holds by part~(ii) of Corollary~\ref{cro:ET}.
On the other hand, we have the inclusion $\Gp_0 \subseteq {}^\perp \proj_n$ by Theorem~\ref{thm:GPT}  and so
we obtain the equality $\Gp_0 \star \proj_n  = \filt (\Gp_0 \cup \proj_n)$ by Lemma~\ref{lem:UV.filt}. 
Therefore, $\Gp_n = \add \big ( \Gp_0 \star \proj_n \big )$. Furthermore, this equality in conjunction with Proposition~\ref{prop:star.operation}
implies that $\Gp_n$ is contravariantly finite provided that $\proj_n$ and $\Gp_0$ are contravariantly finite.

Part~(ii): It is well-known that contravariant finiteness of $\proj_\infty$ implies
$n:= \fpd (\Lambda) < + \infty$; cf.~\ref{prop:artinian.double.perp} or~\cite{AR.applications}.
Thus, $\proj_\infty = \proj_n$ and $\Gp_\infty = \Gp_n$ by Eq.~\eqref{holm.fpd}.
Now by part~(i) of the theorem we have the equality $\Gp_\infty = \add \big ( \Gp_0 \star \proj_\infty \big )$,
which in conjunction with Proposition~\ref{prop:star.operation} implies contravariant finiteness of $\Gp_\infty$.
\eproof

\brem \label{rem:vGor}
It follows from part~(ii) of Theorem~\ref{Main.thmII} that 
for Artin algebras over which the class $\Gp_0$ is contravariantly finite, the following conditions are equivalent:
\begin{itemize}
\item contravariant finiteness of $\proj_\infty$;

\item contravariant finiteness of $\Gp_\infty$.
\end{itemize}
Typical examples of Artin algebras over which the class $\Gp_0$ is contravariantly finite are the so-called \emph{CM-finite algebras} and
\emph{virtually Gorenstein Artin algebras}~\cite{beligiannis.CM,BR}. By~\cite{beligiannis.CM} the latter class of algebras include algebras which are derived equivalent, or stably equivalent of Morita type, with Artin algebras of finite representation type or Gorenstein Artin algebras.
\erem

Finally, we shift our focus from Artin algebras to the slightly more general setting of left artinian rings and prove in the following theorem  that
in parallel to~\cite{trlifaj.approx.little}, contravariant finiteness of $\Gp_\infty$ is still a sufficient condition for validity of the second finitistic dimension conjecture.

\bthm \label{Main.thmIII}
Let $\Lambda$ be a left artinian ring and $(\mathcal{A} , \mathcal{B})$ be the cotorsion pair generated by $\Gp_\infty$.
Let $\{S_1 , \ldots , S_n \}$ be a complete set of simple $\Lambda$-modules and for any $1 \leq i \leq n$ pick a special $\mathcal{A}$-precover $A_i \arrow S_i$.
Then:
\begin{enumerate}
\item $\fpd (\Lambda) = \sup \{ \Gpd_\Lambda (A_i) \mid i = 1 , \ldots , n \}$.

\item $\Gp_\infty$ is contravariantly finite if and only if the modules $A_i$ can be taken finitely generated for every $1 \leq i \leq n$.
In this case, $\fpd (\Lambda) < + \infty$.
\end{enumerate}
\ethm

\bproof

Part~(i): If $m:= \sup \{ \Gpd_\Lambda (A_i) \mid 1 \leq i \leq n \}<+\infty$, then 
$\Gp_\infty \subseteq \mathcal{A} \subseteq \GP_m$ by Proposition~\ref{prop:artinian.double.perp} and the fact that $\GP_m$ is closed under filtrations~\ref{rem:GP}.
Therefore, $\fpd (\Lambda) \leq m$ by Eq.\eqref{holm.fpd}.
If, on the other hand, $m:=\fpd(\Lambda) < + \infty$, then $\Gp_\infty = \Gp_m$ by Eq.~\eqref{holm.fpd} and since $\GP_m$ is closed under filtrations~\ref{rem:GP} it follows that  $\mathcal{A}\subseteq \GP_m$.
Therefore, $\sup \{ \Gpd_\Lambda (A_i) \mid 1 \leq i \leq n \} \leq m$.

Part~(ii): If $\Gp_\infty$ is contravariantly finite, then each $S_i$ has special $\Gp_\infty$-precover $f: A_i \arrow S_i$ in $\mod (\Lambda)$ by~\ref{prop:approx.in.mod}
and so $\ker (f) \in (\Gp_\infty)^\perp \cap \mod (\Lambda) \subseteq \mathcal{B}$. Therefore, $f: A_i \arrow S_i$ is a special $\mathcal{A}$-precover
with $A_i$ finitely generated.
Conversely, if each $A_i$ is finitely generated, then it follows from Corollary~\ref{cro:ET} that $\mathcal{A}^\fin = \Gp_\infty$. In particular each $A_i$ belongs to $\Gp_\infty$ and then it follows from Proposition~\ref{prop:artinian.double.perp} that $\Gp_\infty$ is contravariantly finite.
\eproof

\section*{Acknowledgements}

Major part of this work was done during the first-named author's visit to Department of Algebra at Charles University (MFF UK) in 2018. He wishes to express his gratitude to MFF UK for their hospitality and University of Tehran for financial support of the visit.


\end{document}